\documentclass[Afour,sageh,times]{sagej}

\usepackage{moreverb,url}

\usepackage[colorlinks,bookmarksopen,bookmarksnumbered,citecolor=red,urlcolor=red]{hyperref}

\usepackage{graphicx}
\usepackage{amsmath}
\usepackage{color}
\usepackage{amsthm}
\usepackage{amssymb}

\newcommand\BibTeX{{\rmfamily B\kern-.05em \textsc{i\kern-.025em b}\kern-.08em
T\kern-.1667em\lower.7ex\hbox{E}\kern-.125emX}}

\def\volumeyear{2018}

\begin{document}

\runninghead{Stokowski, Jacak and Jacak}

\title{Effective information processing with pure braid group formalism in view of 2D holographic principle for information}

\author{Hubert Stokowski, Janusz E. Jacak and Lucjan Jacak}

\affiliation{Department of Quantum Technologies, Faculty of Fundamental Problems of Technology, Wroc{\l}aw University of Science and Technology, Wyb. Wyspia\'nskiego 27, 50-370 Wroc{\l}aw, Poland}

\corrauth{Janusz E. Jacak}

\email{janusz.jacak@pwr.edu.pl}

\begin{abstract}
The so-called holographic principle, originally addressed to the high energy physics, suggests more generally that the information inseparatly bounded with a physical carrier (measured by its entropy) scales as the event horizon surface---a two-dimensional object. In this paper we present an idea of representing classical information in formalism of pure braid groups, characterized by an exceptionally rich structure for two-dimensional spaces. This leads to some interesting properties, e.g. information geometrization, multi-character alphabets energetic efficiency for information coding characteristic or a group structure decoding scheme. We also proved that proposed pure braid group approach meets all the conditions for storing and processing of information.
\end{abstract}

\keywords{pure braid group, holographic principle, information processing, information geometrization}

\maketitle

\section{Introduction}

\begin{figure}[!ht]
\centering
\scalebox{0.7}{\includegraphics{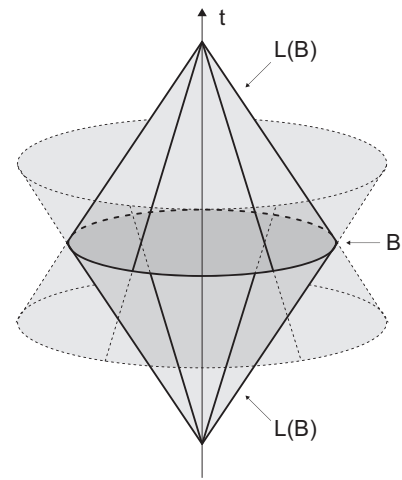}}
\caption{Two light-sheets for 2D spherical surface $B$ model ($B$ is the circle for 2D): $t$ is the time, the cones created by light rays with negative expansion define light-sheets suitable for covariant entropy bound formulated by Bousso: the entropy on each light-sheet $S(L(B))$ does not exceed the area of $B$, i.e., $A(B)$, $S(L(B)) \leq A(B)/4$.}
\label{f1}
\end{figure}

The holographic principle, formulated on the basis of the Bekenstein entropy bound \cite{Bekenstein1972,tHooft2009,Bousso1999}, strongly limits locality \cite{tHooft2009,Bousso1999}. This principle can be formulated \cite{Bousso2002} in terms of a number of independent quantum states describing the light-sheets $L(B)$ (cf. Fig. \ref{f1})---the number of states, ${\cal{N}}$, is bound by the exponential function   of the area $A(B)$ of the surface $B$ corresponding to the light-sheet $L(B)$, 
\begin{equation}
{\cal{N}}[L(B)] \leq e^{A(B)/4}. 
\end{equation}
It can be also expressed equivalently \cite{Bousso2002}, that the number of degrees of freedom (or number of bits multiplied by $\ln 2$ ) involved in description of $L(B)$, cannot exceed $A(B)/4$. Even though the complete holographic theory is not constructed as of yet, this idea is considered as breakthrough on the way to quantum gravity---the holographic principle highlights the lower two-dimensional character of  the information (entropy) corresponding to any real 3D system and in that sense refers to some mysterious 2D hologram collecting all information needed to describe corresponding 3D system. 

An interesting step toward an explicit formulation of the holographic theory was done by Maldacena \cite{Maldacena1997,Witten1998}, who proposed a holographic superstring theory formulation in the anti de Sitter space. According to that model, the 3D quantum gravity could be implemented by 2D quantum nonrelativistic hologram on the edge of the anti de Sitter space (for the anti de Sitter space the boundary is a sphere at fixed time). The locally two dimensional spatial manifolds for the quantum hologram corresponding to 3D quantum relativistic systems open, however, quite new possibility for interpretation of the nature of hypothetical holographic counterparts of elementary particles.

In the context of the holographic principle we propose to grasp this two-dimensional information character with a mathematical formalism favoring two-dimensional spaces over higher-dimensional ones. Such a opportunity provides the braid group formalism, which for two-dimensional spaces is of rich infinite structure unlike for three- or higher-dimensional spaces. This will lead to formulation of some fundamental findings concerning the effectiveness of information coding.

From the informational application point of view we propose to consider a more efficient system for coding information (based on pure braid group formalism) than currently widespread uniformed binary system, which had earned its popular position because of the simplicity of its realization with transistor-based processors. Nowadays this type of information coding is so popular that no other systems are considered to realize in the information processing despite that no type of a binary-code exist in nature.

In this paper we present considerations of multi-character alphabets in coding of classical information, which can be further referred to the realization of neuron-like system of data storage and processing. It seems to be a promising way to enhance the efficiency of the information transmission rate.
	
The main idea of the multi-character alphabets is based on a topological concept of pure braid groups, which has been used in many applications \cite{phd}. Topological properties are highly useful to prove complex theorems in a simpler manner and solve complicated physical problems easier. The original concept of information coding in this field seems to be promising. The idea of storing information with braid groups change the outlook on coding, because it refers not to the commonly known idea of making links one-to-one, which represent a specific information. Braid groups formalism proposed here lets to code information by choosing how do the elements connect instead of which of them are connected. This formalism lets to gain even infinite amount of data stored in a set of $N$ elements, when there are only $N!$ possible bits stored in this type of system in the binary approach. 	

\section{Braid groups in information coding}

\begin{figure}[ht]\centering
\includegraphics[width=4 cm]{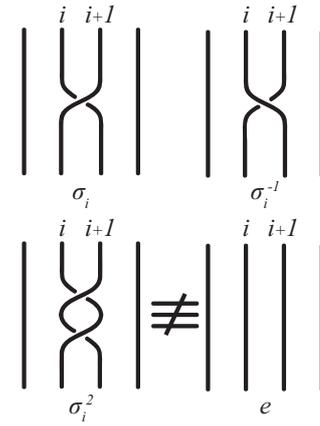} 
\caption{The geometric braid presentation of a full braid group generator, its inverse element and a specific property in the two-dimensional space.}
\label{fig:fig1}
\end{figure}

For 2D plane $R^2$ the braid group was originally described by Artin \cite{Artin1947}, and later vastly developed \cite{Birman2, Birman1, birman3, birman4, Kang, Dehornoy, Guaschi}. The Artin braid group, i.e., $\pi_1(Q_N(R^2))$ is the infinite group, with generators $\sigma_i$ and its inverse $\sigma^{-1}_{i}$ (defining exchanges of neighboring particles, $i$th particle with $(i+1)$th one, at some enumeration of particles, arbitrary, however, due to particle indistinguishability) cf. Fig. \ref{fig:fig1},
\begin{equation}
\sigma_{i}\sigma_{i}\thinspace ^{-1}=e,
\label{eq:equation1}
\end{equation}
satisfying the conditions \cite{Artin1947,Birman1},
\begin{equation}
\sigma_{i} \sigma_{i+1} \sigma_{i} = \sigma_{i+1}\sigma_{i} \sigma_{i+1}, \quad \text{for} \quad 1 \leqslant i\leqslant N-2,
\label{eq:equation2}
\end{equation}
\begin{equation}
\sigma_{i}\sigma_{j}=\sigma_{j}\sigma_{i},	\quad \text{for} \quad 1\leqslant i, j\leqslant N-1, |i-j|\geqslant 2.
\label{eq:equation3}
\end{equation}

Among the different dimensional Euclidean spaces the most valuable is the two-dimensional space. For $dim (M)>2$ the braid group is equal to the permutation group $S_{N}$. The reason is all the closed loops in this type of space can be contracted into a point, no other 0-dimensional topological defect(point) of $M$ of the space can change it. In the 2-dimensional space the loops are uncontractible, thus the braid group is a non-trivial infinite group. This implies different relationships, depending on the dimension of the manifold:
\begin{equation}
(\sigma_{i})^{2}=e, \quad \text{for} \quad \operatorname{dim}(M)>2,
\label{eq:equation4}
\end{equation}
\begin{equation}
(\sigma_{i})^{2}\neq e,  \quad \text{for} \quad \operatorname{dim}(M)=2.
\label{eq:equation5}
\end{equation}

The pure braid group is the first group of homotopy $\pi_{1}(F_{N}(M))=\pi_{1}(M^{N}\setminus\Delta)$, where $F_{N}(M)$ is a configuration space of $N$-braids defined as $n$-fold Cartesian product of manifold $M$ reduced by $\Delta$-points, which represent diagonal points---topological defects of the manifold $M$. In a formal way generators of this group have to be introduced by generators of full braid group $\pi_{1}(Q_{N}(M))=\pi_{1}((M^{N}\setminus\Delta)/S_{N})$, where $S_{N}$ represents a permutation group, which determines the equivalence class in the quotient structure. It means, that the principal difference between full and pure braid groups is that in the full group trajectories can generate braids in any way, with permutations between the beginning and the end, when the pure group is limited by the same setup of trajectories at the beginning and at the end of the braid.

\begin{figure}[ht]\centering
\includegraphics[width=8 cm]{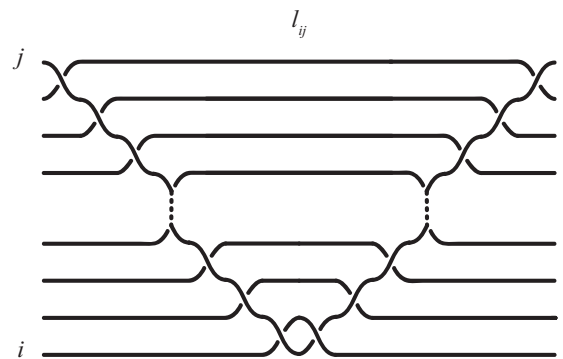} 
\caption{The geometric braid presentation of a pure braid group generator $l_{ij}$.}
\label{fig:fig2}
\end{figure}

The pure braid group is generated by the $l_{ij}$ generators (cf. Fig. \ref{fig:fig2}) formally expressed using elements of the full braid group as:
\begin{equation}
l_{ij}=\sigma_{j-1}\sigma_{j-2}\ldots\sigma_{i+1}\sigma_{i}^{2}\sigma_{i-1}^{-1}\ldots \sigma_{j-2}^{-1}\sigma_{j-1}^{-1}.
\label{eq:equation6}
\end{equation}

Through pure braid groups it seems possible to represent the multi-character alphabet coding system, where each symbol is represented by a different generator of a pure braid group. Representation of this type is useful only in the two-dimensional space, where the structure is non-trivial and the braid is able to contain even infinite countable number of code signs. To proof usability of this coding system it is necessary to show two principal conditions. First of all it is crucial to find at least two (in binary code) or more generators, which are able to represent alphabet characters of the particular alphabet. The second condition is the necessity to preserve the sequence of a code. It has to be shown that there is no possibility to change the order of used generators and it is impossible to untie the braid built by those.

\section{Binary code information}

The basic coding method used nowadays is the binary system. It is principal to show possibility of a building binary system in the formalism of braid groups to proceed with the multi-character alphabets. There is no possibility to show this type of code with the 2-braid, because there is only one generator of the group in this case and its inverse generator:
\begin{equation}
l_{12}=\sigma_{1}^{2}\quad \text{and} \quad l_{12}^{-1}=(\sigma_{1}^{-1})^{2}.
\end{equation}
Building the binary formalism with those two elements is impossible, because the construction of a code with elements inverse to each other would cause disentangle of the braid. The required procedure is adding trajectories to the braid to gain more generators, which can possibly work as characters in the coding alphabets. It can be proven that there exist a possibility to code in the binary system using a 3-braid with a pure braid group. In this type of structure there are 3 possible generators:
\begin{equation}
l_{12}=\sigma_{1}^{2},\quad
l_{13}=\sigma_{2}\sigma_{1}^{2}\sigma_{2}^{-1},\quad
l_{23}=\sigma_{2}^{2}.
\end{equation}
The inverse elements cannot be used to code information in combination with basic generators, thus they are not included here.

Lets prove it is possible to code a binary information through two operators of the 3-braid group. Its realization will be shown by taking two generators, which exchange braids from the group with the specified one. It is trivial in the 3-braid group, nonetheless this proof will be afterwards expanded into more complex structures. Lets take the 3-rd trajectory as the distinguished one and build a segment of generators containing this trajectory. The $l_{13}$ and $l_{23}$ generators can be chosen as possible characters, it has to be proven there is no possibility to change the sequence of the aforementioned just through the braid group properties. In the geometrical braid representation it means that there is no possibility to 'drag' one entanglement through another in the order.

\begin{figure}[ht]\centering
\includegraphics[width=8 cm]{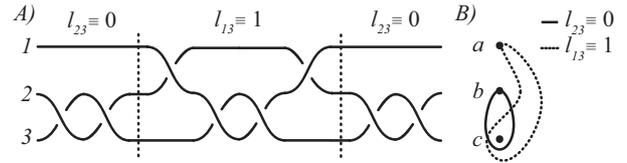} 
\caption{A) The binary code representation through a 3-braid pure group, B) the idea of topological simplification of the problem.}
\label{fig:fig3}
\end{figure}
Showing that the braid generated by sequence of $l_{23}l_{13}l_{23}$ preserves the original input order is enough for verification that any sequence built of those generators do preserve the information structure. The sequence realized by geometrical braids has been shown in the Figure \ref{fig:fig3} A). The proof can be done in two steps. The first one is a topological simplification of the problem, by bringing each unperturbed trajectory into a point and representing each entanglement by a closed loop. In this representation the distinguished, 3-rd trajectory will be transformed into a point and each element of the sequence will be converted into a loop around this point, with the beginning and the end in the point, which represents the second entangled trajectory. It is necessary to show all loops generated by sequence elements are not homotopic and have to cross each other in at least two points. The procedure is shown in the Figure \ref{fig:fig3} B).

\begin{proof}
Suppose $A$ and $B$ are open neighborhoods of the point $c$, which boundaries represents distinguished trajectories and let $a\in \partial A$ but $a \notin cl B$ and analogous $b\in \partial B$ but $b \notin cl A$, where $\partial A$ and $\partial B$ are boundaries of sets $A$ and $B$ respectively. It is therefore true that $c\in A\cap B\neq \emptyset$ and neither $A\nsubseteq B$, nor $B\nsubseteq A$. Both $A$ and $B$ are simply-connected, not border and dense-in-itself sets.

Since $A$ and $B$ are open sets, their intersection $A\cap B$ is open as well, lets take a point $q$, which is a boundary point of $A\cap B$, let $V$ be any open neighborhood of the \textit{q}, then:
\begin{equation}
\forall_{q\in \partial(A\cap B)} \exists_{q', q''\in V} [q'\in (A\cap B) \land  q''\in (X\setminus A \cap B)],
\end{equation}
moreover, due to $\partial (A\cap B)\in \partial A \cup \partial B$:
\begin{equation}
\forall_{q\in \partial(A\cap B)} [q\in \partial A \lor q\in \partial B].
\end{equation}
For all open neighborhoods $V$ of the point $q$ it can be written:
\begin{equation}
\begin{split}
\forall_{q\in \partial(A\cap B)} \forall_{V}  \exists_{q' ,q'' \in V} \{[q'\in A \land q''\in (X\setminus A)]\\
 \lor [q'\in B \land q''\in (X\setminus B)]\}.
 \end{split}
\end{equation}
Since $A$ and $B$ are dense-in-itself and simply-connected, $A\nsubseteq B$, $B\nsubseteq A$ and $A\cap B\neq \emptyset$ it has to be truth, that there exists a boundary point $q_{i}$ of $A\cap B$ and its any open neighborhood $V_{i}$ as:
\begin{equation}
\begin{split}
\exists_{q_{i}} \forall_{i\in (1; 2n)} \forall_{V_{i}} \exists_{q'_{i}, q''_{i}\in V_{i}} [q'_{i}\in (A\cap B) \land q'_{i}\in A \\
\land q'_{i}\in B \land q''_{i}\in X\setminus (A\cap B) \land q''_{i}\in X\setminus A \land q''_{i}\in X\setminus B],
\end{split}
\end{equation}
where $n\in \mathbb{N}$ thus:
\begin{equation}
\exists_{q_{i}} \forall_{i\in (1; 2n)} [q_{i}\in \partial (A\cap B) \land q_{i}\in \partial A \land q_{i}\in \partial B].
\end{equation}
Finally:
\begin{equation}
\partial (A\cap B) \cap \partial A \cap \partial B \neq \emptyset.
\end{equation}
And there are even number of $\partial A$ and $\partial B$ intersections.

It can be as well proven in language of topology and homotopy, what is sketched below. Consider two loops $A$ and $B$ and three points of the manifold $M$: $a$, $b$ and $c$. If the loop $A$ contains the point $a$ and surrounds the point $c$ it cannot be contracted into $a$, due to a defect-like character of the point $c$ in this context. Lets take a second loop $B$ in the analogous manner, but containing the point $b$ instead of $a$. Neither the loop $A$ surrounds the point $b$, nor the loop $B$ surrounds the point $a$ and none of them can be contracted because of the point $c$, which represents a topological defect of the manifold $M$. Moreover none of those loops can be expanded through points $a$ and $b$ for loops $B$ and $A$ respectively. It is clear then, that those loops have to cross each other on the manifold $M$.

\end{proof}

\begin{figure}[ht]\centering
\includegraphics[width= 8cm]{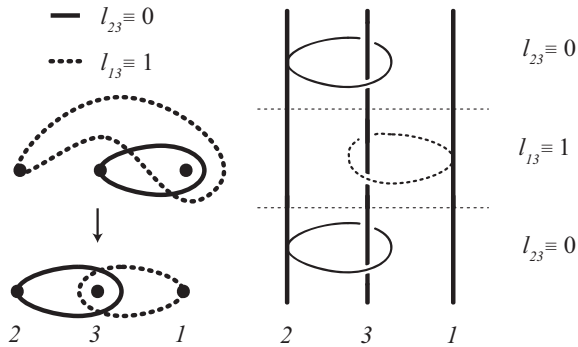} 
\caption{The 3-dimensional representation of the independence of each alphabet symbol of the binary code in the braid groups formalism.}
\label{fig:fig4}
\end{figure}

It is fulfilled that the boundaries representing loops in this approach intersect each other in each case, are not homotopic and cannot be contracted into a point. Basing on the topological proof it is impossible to drag one braid group generator through another in the binary information sequence. In more intuitive way one can undergo with points into the 'bars' in the 3-dimensional space as shown in the Figure \ref{fig:fig4}. It is clear that loops around the bars, which represent different information symbols cannot be mixed with each other.

The other approach to the proof, basing on the braid groups formalism itself can be necessary as well. It has to be shown that there is no possibility to change the order of full braid group generators $\sigma_{i}$, which are used in a realization of pure braid group generators $l_{ij}$ assigned to alphabet characters. As it was said earlier it is enough to show there is no possibility to change the order in  the sequence $l_{23}l_{13}l_{23}$. This sequence can be realized as follows:
\begin{equation}
l_{23}l_{13}l_{23}=\sigma_{2}^{2}\sigma_{2}\sigma_{1}^{2}\sigma_{2}^{-1}\sigma_{2}^{2}.
\end{equation}

Reducing the relation using the neutral element of the group from the equation (\ref{eq:equation1}) and the relation (\ref{eq:equation2}) turns out impossible due to the relation (\ref{eq:equation5}) taking place in the two-dimensional space. Usage of defined relations leads only to complicating the braid, instead of reduction or changing the order of the coded sequence. In result of this paragraph it has been formally shown that there exists a possibility to code information through 3-braid group generators. In this formalism it is possible to use even three generators as independent signs: $l_{13}$, $l_{23}$ and $l_{12}$ to define a ternary system, however the idea described here is easily expandable into more complex braids. The proposed set of generators and assigned alphabet characters are listed in the Table \ref{tab:table1}.

\begin{table}[hbt]\centering
\caption{The binary system defined with 3-braid group generators.}
\begin{tabular}{|c| |c|}
\hline
Braid group generator &
 Binary symbol\\ \hline \hline
$l_{23}$ & 0 \\ 
$l_{13}$ & 1 \\ \hline
\end{tabular}
\label{tab:table1}
\end{table}

\subsection{Ternary system}

\begin{figure}[ht]\centering
\includegraphics[width= 8cm]{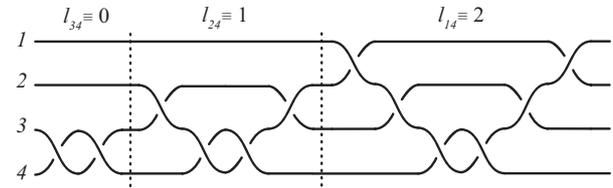} 
\caption{The 4-braid group concept of a ternary coding system.}
\label{fig:fig5}
\end{figure}

\begin{figure}[ht]\centering
\includegraphics[width= 8cm]{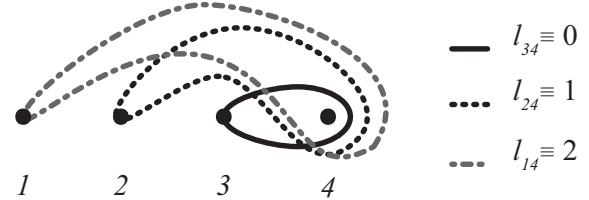} 
\caption{The scheme for a topological proof for the independence of ternary system generators in a 4-braid system.}
\label{fig:fig6}
\end{figure}

Expanding the proposed approach it is easy to define a ternary system by adding one more trajectory to the braid. Lets define a 4-braid with 3 independent generators, which can be used as symbols. In the new 4-braid system by distinguishing the 4-th trajectory one can choose at least 3 independent generators, the concept is shown in the Figure \ref{fig:fig5}.
The presented above proof can be easily expanded into the 4-braid situation and show the independence of 3 mentioned generators. The only difference for systems with number of trajectories $N\geqslant 4$ is that one has to put attention to relation (\ref{eq:equation3}). Its consequence is impossibility to choose two generators tangling up completely different trajectories. Choosing this type of generators would cause a failure to preservation of the information sequence - dragging those braids would be unrestricted and will lead to wasting of the sequence order. However, this disadvantage  does not affect the proposed coding system. The generators of the chosen ternary system has been shown in the Table \ref{tab:tab2}.

\begin{table}[hbt]\centering
\caption{The ternary system defined with 4-braid group generators.}
\begin{tabular}{|c| |c|}
\hline
Braid group generator &
 Ternary symbol\\ \hline \hline
$l_{34}$ & 0 \\ 
$l_{24}$ & 1 \\ 
$l_{14}$ & 2 \\ \hline
\end{tabular}
\label{tab:tab2}
\end{table}
For an intuitive topological proof it is comfortable to base on the scheme shown in the Figure \ref{fig:fig6} and expand the proof for 3 neighborhoods of a point-representation of the distinguished trajectory. For a braid group formalism it is enough to show the sequence of $l_{34}l_{24}l_{14}l_{24}l_{34}l_{14}l_{34}$ preserves the order or just to show that each pair of the generators ($l_{34}l_{24}$, $l_{24}l_{34}$, $l_{24}l_{14}$, $l_{14}l_{24}$ etc.) does preserve it.

\subsection{Multi-character alphabets}

As it has been sketched above in introduced formalism it is possible to expand coding alphabets by adding more and more trajectories to the braid and each addition has to expand amount of generators eligible for coding information at least by one. The consequence is very strong and leads to the proposal of the multi-character alphabet definition using the 2-dimensional braid group formalism, where number of alphabet signs increases at least linearly with the number $N$ of trajectories in the \textit{N}-braid. The discussion of adding more than one generator with each new trajectory is open, but the main goal of this work is to show there is a possibility to expand multi-character alphabets in the countable and infinite manner. Adding more generators is possible, as it has been shown in the 3-braid case, but it is hard to generalize into more complex braid structures, mainly due to the relation (\ref{eq:equation3}). Possible representation of the \textit{N}-character alphabet coded using the \textit{(N+1)}-braid group with the \textit{(N+1)}-th trajectory distinguished has been shown in the Table \ref{tab:tab3}.

\begin{table}[hbt]\centering
\caption{The \textit{N}-character alphabet coded using the \textit{(N+1)}-braid group formalism.}
\begin{tabular}{|c| |c|}
\hline
Braid group generator &
 \textit{N}-alphabet symbol\\ \hline \hline
$l_{N N+1}$ & 0 \\ 
$l_{N-1 N+1}$ & 1 \\ 
$l_{N-2 N+1}$ & 2 \\ 
\vdots & \vdots \\ 
$l_{1 N+1}$ & N-1 \\ \hline
\end{tabular}
\label{tab:tab3}
\end{table}

\section{Metrics and decoding distance}

\begin{figure}[ht]\centering
\includegraphics[width= 8cm]{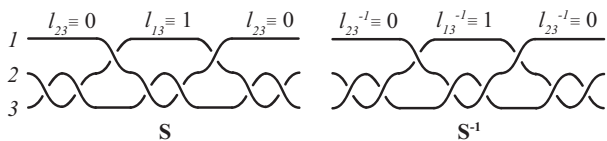} 
\caption{The pure braid group information decoding idea in terms of geometric braids.}
\label{fig:fig7}
\end{figure}

In the perspective of applications the decoding case can be considered. The most intuitive way in the pure braid groups formalism, as it has a group structure, is to use the full string of inverse generators. Let $S$ be a string of \textit{d} binary-symbols in a pure braid sub-group representation and $S^{-1}$ be its inverse string as shown in the Figure \ref{fig:fig7}.

Through acting on the string $S$ with its inverse $S^{-1}$ it can be found that the whole constructed braid is an equivalent to the neutral element of the group $e$---decoded string has disentangled completely with its inverse:
\begin{equation}
S\cdot S^{-1}\equiv l_{23}l_{13}l_{23}\cdot l_{23}^{-1}l_{13}^{-1}l_{23}^{-1}\equiv e.
\end{equation}
In general the string $S$ of a coded information can be represented in a formal way with pure braid group generators:
\begin{equation}
S={\displaystyle \prod_{j\in a_{d}}^{} l_{jn}},
\end{equation}
where $a_{d}$ is a sequence build of natural numbers corresponding to specific alphabet signs. The product $\prod_{j\in a_{d}}^{} l_{jn}$ is understood as a combination of following pure braid group generators. The inverse string $S^{-1}$ can be defined as:
\begin{equation}
S^{-1}={\displaystyle \prod_{j\in a'_{d}}^{} l_{jn}^{-1}},
\end{equation}
where the sequence $a'_{d}$ is given by:
\begin{equation}
a'_{d}:\, \forall_{i\in (1,d)}\, a'_{i}=a_{(d+1-i)}.
\end{equation}

In terms of the information theory it is common to use specific types of metrics, so called information distance. Using the classical binary system there is a well defined Hamming distance \cite{Hamming}, which is constructed between two equal strings as a number of positions at which the corresponding symbols are different. Its generalization on any string lengths has been shown by Levenstein \cite{Levenstein}. In the case of braid-groups it is more convenient to define a different metrics, corresponding to the proposed decoding method.

Lets consider a metrics, which correspond to the number of equal characters at the end of two \textit{d}-character strings built according to the specific braid sub-group with the distinguished \textit{n}-th trajectory. In case of two strings, which lengths are not equal lets add its difference to the resultant distance between information coded on those. The metrics corresponding to the inverse-braid decoding system for two strings $S_{1}$ and $S_{2}$ can be proposed as:
\begin{equation}
d(S_{1},S_{2})=\cfrac{f+d+|f-d|-2x}{2},
\label{eq:equation7}
\end{equation}
where $d$ is the string's $S_{1}$ length, $f$ is the string's $S_{2}$ length and $x$ is a number of repeating characters at the end of those, formally:
\begin{equation}
S_{1}={\displaystyle \prod_{j\in a_{d}}^{} l_{jn}}, \qquad
S_{2}={\displaystyle \prod_{k\in b_{f}}^{} l_{kn}},
\end{equation}
$a_{d}$ and $d_{f}$ are \textit{d}-element and \textit{f}-element sequences consisting of natural numbers, corresponding with code characters coded with $l_{jn}$ and $l_{kn}$ - pure braid group generators acting with the \textit{n}-th distinguished trajectory. Moreover $\xi_{x}$ is a subsequence of both $a_{d}$ and $b_{f}$:

\begin{equation}
x:\quad \xi_{x}:\quad \xi_{x}\in a_{d}\  \land \ \xi_{x}\in b_{f},\ \text{thus}:
\label{eq:sequence}
\end{equation}
\begin{equation}
\begin{split}
{\displaystyle \forall_{i\in (1;x)}\ \xi_{i}=a_{(d-x+i)}\ \land \xi_{i}=b_{(f-x+i)}\implies }\\
{\displaystyle \forall_{i\in (1;x)}\ a_{(d-x+i)}=b_{(f-x+i)}}.
\end{split}
\end{equation}

\begin{proof}
It has to be proven that the proposed metrics meets four basic conditions:
\paragraph{The non-negativity or separation axiom.}

If $\xi_{x}$ is as well a subsequence of $a_{d}$ as $b_{f}$, so $x\leqslant d$  and $x\leqslant f$, thus $f+d\geqslant 2x$ and $d(S_{1},S_{2})\geqslant 0$ is always fulfilled.

\paragraph{The identity of indiscernibles.}
\begin{equation}
d(S_{1},S_{2})=0\Leftrightarrow f+d+|f-d|-2x=0.
\end{equation}
For $f\in \mathbb{N}$ and $d\in \mathbb{N}$, which is determined by sequences $a_{d}$ and $b_{f}$ properties, expression $f+d-|f-d|$ can be interpreted as doubled maximum of a set $(d, f)$.
\begin{equation}
\begin{split}
f+d+|f-d|=2max(d, f), \text{so:}\\
max(d,f)=x, \text{thus:}\\
x=f\Leftrightarrow f\geqslant d \lor x=d\Leftrightarrow d\geqslant f,
\end{split}
\end{equation}
moreover from (\ref{eq:sequence}):
\begin{equation}
x\leqslant f \land x\leqslant d.
\end{equation}
Finally one gets:
\begin{equation}
\begin{split}
x=f=d\implies\\
\forall_{i\in (1;x)} a_{i}=b_{i}\implies a_{d}\equiv b_{f}\implies\\
d(S_{1},S_{2})=0\Leftrightarrow S_{1}=S_{2}.
\end{split}
\end{equation}

\paragraph{Symmetry.}

According to the definition:
\begin{equation}
d(S_{2},S_{1})=\frac{d+f+|d-f|-2x}{2},
\end{equation}
is always equal to:
\begin{equation}
d(S_{1},S_{2})=\frac{f+d+|f-d|-2x}{2}.
\end{equation}
Taking into account assumptions (\ref{eq:sequence}):
\begin{equation}
S_{2}={\displaystyle \prod_{k\in b_{f}}^{} l_{kn}}, \qquad
S_{1}={\displaystyle \prod_{j\in a_{d}}^{} l_{jn}}.
\end{equation}
\begin{equation}
\begin{split}
x:\quad \xi_{x}:\quad \xi_{x}\in b_{f}\  \land \ \xi_{x}\in a_{d},\ \text{thus}:\\
{\displaystyle \forall_{i\in (1;x)}\ \xi_{i}=b_{(f-x+i)}\ \land \xi_{i}=a_{(d-x+i)}\implies }\\
{\displaystyle \forall_{i\in (1;x)}\ b_{(f-x+i)}=a_{(d-x+i)}}\Leftrightarrow\\
{\displaystyle \forall_{i\in (1;x)}\ a_{(d-x+i)}=b_{(f-x+i)}},
\end{split}
\end{equation}
finally:
\begin{equation}
\forall_{(S_{1},S_{2)}} d(S_{1},S_{2})=d(S_{2},S_{1}).
\end{equation}

\paragraph{Triangle inequality.}

To proof the triangle inequality three different strings have to be defined:
\begin{equation}
S_{1}={\displaystyle \prod_{j\in a_{d}}^{} l_{jn}}, \qquad
S_{2}={\displaystyle \prod_{k\in b_{f}}^{} l_{kn}}, \qquad
S_{3}={\displaystyle \prod_{m\in c_{g}}^{} l_{mn}},
\end{equation}
$a_d$, $b_f$ and $c_g$ are three different length sequences of natural numbers, corresponding with alphabet characters. In general each pair of sequences is able to have different number of the same signs at the end, lets take x, y and z.
For $S_{1}$ and $S_{3}$:
\begin{equation}
x:\quad \xi_{x}:\quad \xi_{x}\in a_{d}\  \land \ \xi_{x}\in c_{g},
\label{eq:sequence13}
\end{equation}
\begin{equation}
\begin{split}
{\displaystyle \forall_{i\in (1;x)}\ \xi_{i}=a_{(d-x+i)}\ \land \xi_{i}=c_{(g-x+i)}\implies }\\
{\displaystyle \forall_{i\in (1;x)}\ a_{(d-x+i)}=c_{(g-x+i)}},
\end{split}
\end{equation}
for $S_{1}$ and $S_{2}$:
\begin{equation}
y:\quad \xi'_{y}:\quad \xi'_{y}\in a_{d}\  \land \ \xi'_{y}\in b_{f},
\label{eq:sequence12}
\end{equation}
\begin{equation}
\begin{split}
{\displaystyle \forall_{i\in (1;y)}\ \xi'_{i}=a_{(d-y+i)}\ \land \xi'_{i}=b_{(f-y+i)}\implies }\\
{\displaystyle \forall_{i\in (1;y)}\ a_{(d-y+i)}=b_{(f-y+i)}},
\end{split}
\end{equation}
for $S_{2}$ and $S_{3}$:
\begin{equation}
z:\quad \xi''_{z}:\quad \xi''_{z}\in b_{f}\  \land \ \xi''_{z}\in c_{g},
\label{eq:sequence23}
\end{equation}
\begin{equation}
\begin{split}
{\displaystyle \forall_{i\in (1;z)}\ \xi''_{i}=b_{(f-z+i)}\ \land \xi''_{i}=c_{(g-z+i)}\implies }\\
{\displaystyle \forall_{i\in (1;z)}\ b_{(f-z+i)}=c_{(g-z+i)}}.
\end{split}
\end{equation}
It can be written that:
\begin{equation}
\begin{split}
d(S_{1},S_{2})+d(S_{2},S_{3})=\\
\cfrac{f+d+|f-d|-2y}{2}+\cfrac{g+f+|g-f|-2z}{2}=\\
\cfrac{g+d+|g-f|+|f-d|+2f-2y-2z}{2}.\\
\end{split}
\end{equation}
From the absolute value properties:
\begin{equation}
\begin{split}
\forall_{x_{i}\in \mathbb{R}} |x_{1}+x_{2}|\leqslant |x_{1}|+|x_{2}|,\ \text{for}\ x_{1}=g-f,\; x_{2}=-d+f\\
|g-f|+|f-d|\geqslant |g-d|,\ \text{thus:}\\
\cfrac{g+d+|g-f|+|f-d|+2f-2y-2z}{2}\geqslant\\
\geqslant \cfrac{g+d+|g-d|}{2}+(f-y-z).
\end{split}
\end{equation}
Two cases have to be considered. According to the definition, sequences $a_{d}$ and $c_{g}$ are the same on the \textit{x} closing elements. If $a_{d}$ and $b_{f}$ are the same on the \textit{y} last elements and:

1) $y\leqslant x$ then $z=y\implies z\leqslant x$ has to be fulfilled, because the same amount of closing elements have to be equal in pairs $a_{d}$, $b_{f}$ and $b_{f}$, $c_{g}$. Moreover according to (\ref{eq:sequence12}) $y\leqslant f$, so $(f-y)\geqslant 0$. Then $(f-y-z)\geqslant (-x)$.

2) $y> x$ then $z=x$, to fulfill assumption that only x closing elements of $a_{d}$ and $c_{g}$ are equal. $(f-y)\geqslant 0$ is conserved, so $(f-y-z)=(-x)$.

In both cases:
\begin{equation}
\begin{split}
\cfrac{g+d+|g-d|}{2}+(f-y-z)\geqslant \\
\geqslant \cfrac{g+d+|g-d|}{2}-x=d(S_{1},S_{3}),
\end{split}
\end{equation}
thus:
\begin{equation}
\forall_{(S_{1},S_{2},S_{3})} d(S_{1},S_{3})\leqslant d(S_{1},S_{2})+d(S_{2},S_{3}).
\end{equation}

\end{proof}
With the well defined metrics it is useful to consider a distance between the set of exemplary 3-bit sequence. The comparison has been shown in the Table \ref{tab:tab4}.

\begin{table}[hbt]\centering
\caption{The comparison of information distances between all 3-bit codes.}
\scriptsize
\begin{tabular}{|c| |c| |c| |c| |c| |c| |c| |c|}
\hline
\textbf{000} & \textbf{001} & \textbf{011} & \textbf{100} & \textbf{010} & \textbf{110} & \textbf{101} & \textbf{111} \\ \hline \hline
$\textbf{001}\,\, \textit{1}$ & $\textbf{000}\,\, \textit{1}$ & $\textbf{011}\,\, \textit{1}$ & $\textbf{010}\,\, \textit{1}$ & $\textbf{101}\,\, \textit{1}$ & $\textbf{111}\,\, \textit{1}$ & $\textbf{100}\,\, \textit{1}$ & $\textbf{110}\,\, \textit{1}$ \\ \hline
$\textbf{010}\,\, \textit{2}$ & $\textbf{010}\,\, \textit{2}$ & $\textbf{000}\,\, \textit{2}$ & $\textbf{000}\,\, \textit{2}$ & $\textbf{110}\,\, \textit{2}$ & $\textbf{100}\,\, \textit{2}$ & $\textbf{110}\,\, \textit{2}$ & $\textbf{100}\,\, \textit{2}$ \\ \hline
$\textbf{011}\,\, \textit{2}$ & $\textbf{011}\,\, \textit{2}$ & $\textbf{001}\,\, \textit{2}$ & $\textbf{001}\,\, \textit{2}$ & $\textbf{111}\,\, \textit{2}$ & $\textbf{101}\,\, \textit{2}$ & $\textbf{111}\,\, \textit{2}$ & $\textbf{101}\,\, \textit{2}$ \\ \hline
$\textbf{100}\,\, \textit{3}$ & $\textbf{100}\,\, \textit{3}$ & $\textbf{100}\,\, \textit{3}$ & $\textbf{100}\,\, \textit{3}$ & $\textbf{011}\,\, \textit{3}$ & $\textbf{011}\,\, \textit{3}$ & $\textbf{011}\,\, \textit{3}$ & $\textbf{011}\,\, \textit{3}$ \\ \hline
$\textbf{110}\,\, \textit{3}$ & $\textbf{110}\,\, \textit{3}$ & $\textbf{110}\,\, \textit{3}$ & $\textbf{110}\,\, \textit{3}$ & $\textbf{001}\,\, \textit{3}$ & $\textbf{001}\,\, \textit{3}$ & $\textbf{001}\,\, \textit{3}$ & $\textbf{001}\,\, \textit{3}$ \\ \hline
$\textbf{101}\,\, \textit{3}$ & $\textbf{101}\,\, \textit{3}$ & $\textbf{101}\,\, \textit{3}$ & $\textbf{101}\,\, \textit{3}$ & $\textbf{010}\,\, \textit{3}$ & $\textbf{010}\,\, \textit{3}$ & $\textbf{010}\,\, \textit{3}$ & $\textbf{010}\,\, \textit{3}$ \\ \hline
$\textbf{111}\,\, \textit{3}$ & $\textbf{111}\,\, \textit{3}$ & $\textbf{111}\,\, \textit{3}$ & $\textbf{111}\,\, \textit{3}$ & $\textbf{000}\,\, \textit{3}$ & $\textbf{000}\,\, \textit{3}$ & $\textbf{000}\,\, \textit{3}$ & $\textbf{000}\,\, \textit{3}$ \\ \hline
\end{tabular}
\label{tab:tab4}
\end{table}

\begin{figure}[ht]\centering
\includegraphics[width=5 cm]{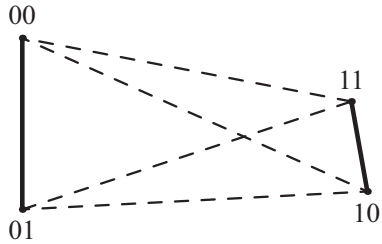} 
\caption{The geometrical representation of distances between all possible 2-bit information using the proposed metrics. Solid line represent the distance of d=1, dashed line - d=2.}
\label{fig:fig8}
\end{figure}

The Figure \ref{fig:fig8} shows a graphical intuition for understanding distances in a simpler case - between all possible 2-bit binary sequences. To represent a 3-bit code it would be necessary to use a four-dimensional picture.

The most valuable conclusion is for every 3-bit code there is only one code in a distance of unity, two with a distance of 2 and 4 with distance of 3. It is the significant difference between proposed metrics and e.g. the Hamming distance, where more codes are 'closer' to each other. It is a valuable conclusion due to the error detection and the error correction issues \cite{Wagner}. In terms of the Hamming distance for d=1 it is neither possible to detect nor correct errors. For d=2 it is possible to detect them and for d=3 also the correction is achievable. If it is possible to use those relations in terms of the proposed decoding system and metrics it can turn out much more effective.

\section{Coding Efficiency}

\begin{figure}[ht]\centering
\includegraphics[width=8 cm]{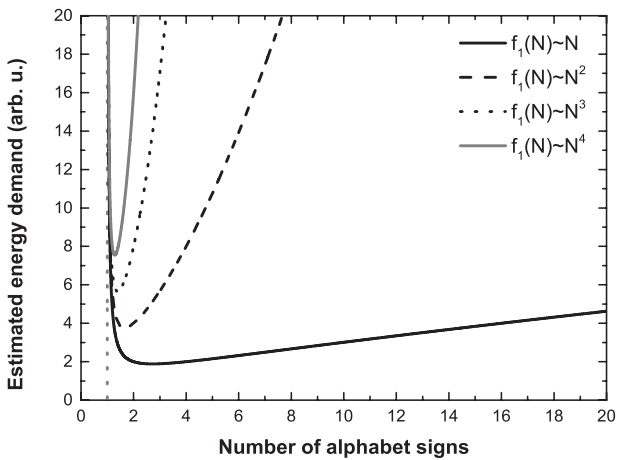} 
\caption{The energy consumption versus number of alphabet signs for $f_{1}(N)$ proportional to $N^{i}$ for $i=1,2,3,4$.}
\label{fig:fig9}
\end{figure}

\begin{figure}[ht]\centering
\includegraphics[width=8 cm]{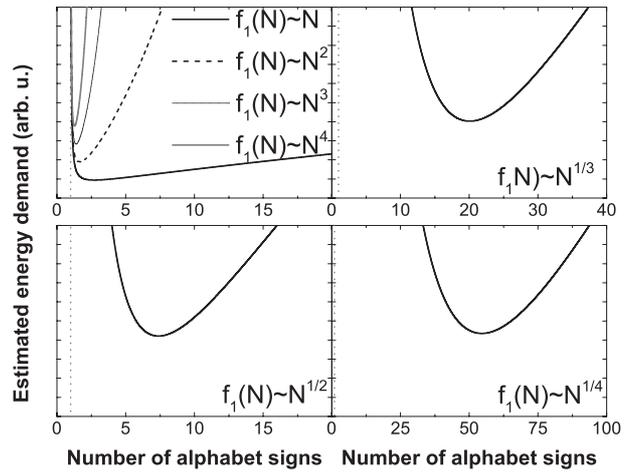} 
\caption{The energy consumption versus number of alphabet signs for $f_{1}(N)$ proportional to $N^{i}$ for $i=1/2, 1/3, 1/4$.}
\label{fig:fig10}
\end{figure}

Coding multi-character alphabets, which could be used to store and process information in the 2-dimensional pure braid group formalism with possibly infinite number of characters, engenders a new issue for the optimization of the number of characters in the coding system. It is important to find an optimal number of alphabet signs for minimization of the energy efforts. This type of optimum should exist, because the greater amount of characters leads to the reduction of the information length and, on the other hand, adding more and more characters requires increasing of the energy consumption to 'remember' all of them. One can try to model an analytic function suitable to predict the energy effort necessary to take an advantage of N-character alphabets versus the coding gain estimated for this amount of available signs, let $g(N)$ be this function. To find an optimum number \textit{N} of alphabet signs it is important to take a function $f_{1}(N)$, which represent the energy efforts for a \textit{N}-character alphabet and divide it by a function $f_{2}(N)$, modeling the gain vs. the binary system efficiency. $f_{1}(N)$ is an increasing function, because using higher number of alphabet signs has to utilize the higher energy consumption. $f_{2}(N)$ can be represented, with respect to the basics of the information theory as:
\begin{equation}
f_{2}(N)=\log_{2}N.
\label{eq:equation12}
\end{equation}
There is no formal function, matching the $f_{1}(N)$ requirements, but if it is indeed increasing the resultant $g(N)$ function has to have one, trivial minimum for $N\geqslant 1$.
\begin{equation}
g(N)=\frac{f_{1}(N)}{f_{2}(N)}.
\label{eq:g(N)}
\end{equation}
Some increasing functions can be proposed to examine the course of the function (\ref{eq:g(N)}), a few basic ideas has been shown on the Figure \ref{fig:fig9}, for $f_{1}(N)=N^{i}$, where $i=1,2,3,4$ and on the Figure \ref{fig:fig10} for $i=\frac{1}{2}, \frac{1}{3}, \frac{1}{4}$. For $f_{1}(N)\propto \sqrt[3]{N}$ the optimum number of alphabet signs \textit{N} is about 20, which naively corresponds to the number of syllables in the most of the languages commonly used by humans. It can be a clue for the naturally preferable amount of characters in the alphabet and also for a possibly accurate model of the energy consumption vs. number of alphabet signs as $f_{1}(N)\propto \sqrt[3]{N}$ in natural systems (cube root corresponds to a 3D volume).

\section{Conclusion}

Presented approach emphasizes that every information must be connected with its physical carrier, and this relation (between information and its carrier) allows to discuss for example energy balance of a information processing system (in case when one considers information without its carrier then in fact the binary system seems to be most optimal due to its simplicity). On the other hand, the carrier must behave accordingly to its physical properties, which in view of a general character of the holographic principle, leads to quite interesting observation that every physical information (information bounded to its carrier---in fact only such information exists) scales as a two-dimensional object and not a three-dimensional one, as it could be falsely suspected. Thus one can suspect that it would be desirable and effective to analyze information and its processing systems with use of some mathematical description which favors two-dimensional spaces among all others.

A multi-character alphabet coding using pure braid groups has been introduced in context of the holographic principle elucidating two-dimensional information character (where complex properties of pure braid groups are only preserved in the $\operatorname{dim}(M)=2$). If properly implemented this could significantly improve the coding efficiency in the information technology. As it is hard to find in nature a truly binary information processing system (other than artificially introduced systems in computer science), a multi-character alphabet based information processing system seems to be more natural. Such system turns out more energy-efficient than the binary system, which is commonly used nowadays. The proposed information decoding system forces an usage of a specific metrics, however it can be an advantage of the described information coding due to its potential in the error-detection and the error-correction issue, since distances between the sequences are usually high. The proposed formalism gives a possibility to connect the information with its physical carrier---the braid, which has not been proposed so far. Moreover, the described coding system seems to be a suitable way for a realization of neuron-like resonance circuits for the data storage and is likely to explain working principles of this type of connections in the human brain.

The method described above introduces just one of possibilities for the coding of information with pure braid groups generators. An equivalent method is using generators responsible for tangling up the distinguished trajectory with other ones. The limit imposed on the pure braid group generators can be as well a defect of this approach, because of the usage of a sub-group, as an advantage. It could be considered to preform an efficient error-correction system using a few equivalent representations of the alphabet transferring through different transmission channels. Other approaches can find more than $N-1$ generators useful for the coding in the $N$-braid, as it has been shown for the 3-braid, but the method proposed here cannot generalize this type of generators due to the additional relation (\ref{eq:equation3}), appearing in systems with $N\geqslant 4$ trajectories. Another way to earn more possible characters for the coding system is a realization of more complicated form of generators, which would be able to preserve the necessary conditions.

\begin{acks}
Supported by the NCN project P.2011/02/A/ST3/00116
\end{acks}

\end{document}